\documentclass[12pt]{article}%
\usepackage{amsmath}
\usepackage{amsfonts}
\usepackage{amssymb}
\usepackage{graphicx}
\usepackage[a4paper]{geometry}%
\newtheorem{theorem}{Theorem}

\newtheorem{corollary}[theorem]{Corollary}

\newtheorem{lemma}[theorem]{Lemma}

\newtheorem{remark}[theorem]{Remark}

\newenvironment{proof}[1][Proof]{\noindent\textbf{#1.} }{\ \rule{0.5em}{0.5em}}
\DeclareMathOperator{\diver}{div}
\DeclareMathOperator{\rotor}{rot}
\DeclareMathOperator{\grad}{grad}
\begin{document}

\title{On duality of spaces of harmonic vector fields}
\author{Ren\'{e} D\'{a}ger\footnote{Departamento de Fundamentos Matem\'aticos, ETSIA, Universidad
Polit\'ecnica de Madrid, Plaza Cardenal Cisneros 3, 28040 Madrid,
Spain (rene.dager@upm.es).} and Arturo Presa\footnote{Mathematics
Department, Miami Dade College, Wolfson Campus, Miami, FL33132-2296,
USA (apresa@mdc.edu).}} \date{} \maketitle

\abstract{A differential form $u$ of class $C^{\infty}$ on the
Riemannian manifold $X$ is said to be harmonic if it is closed and
co-closed, i.e.,
\[
du=0,\qquad\delta u=0.
\]

Harmonic differential forms are a natural multi-dimensional
extension of the concept of analytic function of complex variable.
In this paper we characterize continuous linear functionals acting
of the space of germs of harmonic differential forms on a compact
set. This result provides a multi-dimensional analog of a theorem by
G. K\"othe on the dual of the space of germs of analytic functions
of complex variables on a compact.}
\section{Introduction}

A vector field $\bar{u}$ defined in an open set $\mathcal{O}\subset%
\mathbb{R}
^{3}$ is said to be harmonic if it satisfies the equations%
\[%
\rotor
\bar{u}=0,\qquad%
\diver
\bar{u}=0
\]
in $\mathcal{O}$.

The analogy between harmonic vector fields and analytic functions is
explained by the fact that both objects are particular cases of the
concept of harmonic differential form (see Introduction in
\cite{pre}). In this work a differential form $u$ of class
$C^{\infty}$ on the Riemannian manifold $X$ is said to harmonic if
it is closed and co-closed, i.e.,
\[
du=0,\qquad\delta u=0.
\]

It was proved by Ko\"{e}the (\cite{koe}) that the dual space of the space
$h(K)$ of germs of analytic functions on a compact set $K\subset%
\mathbb{C}
$ with the inductive limit topology is isomorphic to the space of analytic
functions on the complementary of $K$ that vanish at infinity, that is%
\[
\left(  h(K)\right)  ^{\ast}=h_{0}(K^{c}).
\]

One example of the significance of this fact is the proof given by
Havin of a formula analogous to the Laurent series expansion for
analytic functions vanishing at infinity with singularities
concentrated in a compact set.

Thus, it is natural to study the problem of the description of the
dual of the space of germs of harmonic vector fields. It turns out
that the space that provides the answer to this question is not
precisely the space of harmonic vector fields on $K^{c}$ that vanish
at infinity. Hence a second problem arises: to find a space whose
dual is isomorphic to $h_{0}(K^{c})$.

In this paper both questions  above are studied in a more general
context, namely, for harmonic differential forms on an Euclidean
space of finite dimension. The main results for the case $n=3$ were
announced in \cite{dp}.

The spaces that solve the mentioned problems are spaces of certain
classes of holomorphic pairs, i.e., non-homogeneous differential
forms $w_{r-1}+w_{r+1}$, where $dw_{r-1}+\delta w_{r+1}=0,\delta
w_{r-1}=0,dw_{r+1}=0$. Such forms where introduced by Gustaffson and
Havinson \cite{gus} in connection with the characterization of
annihilators of harmonic forms on a smoothly bounded domain. For
$n=3$ and $r=1$ holomorphic pairs coincide with holomorphic vector
fields.

As an application of our results we prove in Section 5.B that every
harmonic form defined on the complementary of a compact set may be
represented as the sum of a differential and a codifferential of
suitable forms.

\section{Preliminaries}

\subsection{Some spaces of differential forms and linear functionals}

Throughout this paper $E$ denotes an oriented Euclidean space of dimension
$n\geq2$. The symbol $%
{\textstyle\bigwedge^{r}}
E$  denotes the space of $r$-covectors on $E$. Let also for $r=0,$ $%
{\textstyle\bigwedge\nolimits^{0}}
E=%
\mathbb{R}
$ and $$%
{\textstyle\bigwedge^{\ast}}
E=%
{\textstyle\bigoplus_{k=0}^{n}} {\textstyle\bigwedge^{r}} E.$$

In $%
{\textstyle\bigwedge^{\ast}}
E$ an exterior product is introduced that turns this space into an
algebra. The main properties of this algebra may be found in
\cite{fed}. The $\ast$ operator (the Hodge conjugate) acting on
$r$-covectors is understood in the sense of \cite{fed}. In that
reference  a suitable notion of scalar product of $r$-covectors may
be also found.

\textbf{A. }Here we gather  some facts on  harmonic differential
forms. First, recall the definition of the
operator $\delta$:%
\[
\delta u:=(-1)^{nr+r+1}\ast d\ast u,
\]
where $u$ is an $r$-form of class $C^{1}$, $\ast$ is the Hodge conjugate
operator and $d$ is the differential of forms.

Below we define some spaces of differential forms that we shall use in this paper.

Let $H_{r}(\mathcal{O})$ denote the space of differential forms $u$ of degree
$r$ and of class $C^{1}$ in $\mathcal{O}$ and such that
\[
du=0,\qquad\delta u=0
\]
in $\mathcal{O}$. The elements of $H_{r}(\mathcal{O})$ are called
\emph{harmonic differential forms of degree }$r$\emph{ defined in}
$\mathcal{O}$.

Let $Wh_{r}(\mathcal{O})$ be the space of differential forms of degree $r$ and
of class $C^{2}$ in $\mathcal{O}$ such that
\[
\Delta u:=-(d\delta+\delta d)u=0
\]
in $\mathcal{O}$. The forms belonging to $Wh_{r}(\mathcal{O})$ are called
\emph{weakly harmonic forms defined in} $\mathcal{O}$.

The space $P_{r}(\mathcal{O})$ for $1\leq r\leq n-1$ consists of
non-homogeneous forms $u=w_{r-1}+w_{r+1}$, where $w_{r-1}$ and $w_{r+1}$ are
forms of degree $r-1$ and $r+1$, respectively, which are of class $C^{1}$ in
the open set $\mathcal{O}$ and such that the following equations are valid in
$\mathcal{O}$:%
\[
dw_{r-1}+\delta w_{r+1}=0,\quad dw_{r+1}=0,\quad\delta w_{r-1}=0.
\]
The forms belonging to $P_{r}(\mathcal{O})$ are called \emph{holomorphic pairs
corresponding to }$r$\emph{ defined in} $\mathcal{O}$.

By means of $Wh(\mathcal{O})$ we denote the direct sum $\bigoplus_{r=0}%
^{n}Wh_{r}(\mathcal{O})$. Clearly, $H_{r}(\mathcal{O})\subset Wh_{r}%
(\mathcal{O})$ and $P_{r}(\mathcal{O})\subset Wh(\mathcal{O})$.

In each of the mentioned spaces the uniform convergence topologies are
considered. In $H_{r}(\mathcal{O})$, for instance, this topology is obtained
from the system of neighborhoods of null%
\[
\left\{  V(K,\varepsilon):\text{ }K\text{ compact and }K\subset\mathcal{O}%
,\varepsilon>0\right\}
\]
where $V(K,\varepsilon):=\left\{  u\in H_{r}(\mathcal{O}):\max_{K}\left\Vert
u(x)\right\Vert <\varepsilon\right\}  $. Here $\left\Vert u\right\Vert
=\sqrt{\left\langle u,u\right\rangle }$ and $\left\langle \cdot,\cdot
\right\rangle $ is the scalar product in $%
{\textstyle\bigwedge^{r}}
E$. The same happens with the spaces $Wh_{r}(\mathcal{O})$ and $P_{r}%
(\mathcal{O})$.

Let $H_{r,0}(E\setminus K)$ denotes, for a compact set $K$, the subspace of
$H_{r}(E\setminus K)$ whose elements vanish at infinity. In the same way it is
defined $P_{r,0}(E\setminus K)$ as a subspace of $P_{r}(E\setminus K)$.

Moreover, $H_{r}(K),Wh_{r}(K),Wh(K)$ and $P_{r}(K)$ represent the
corresponding inductive limits spaces. In order to avoid repetitions let us
see more precisely just the definition of $H_{r}(K)$. The elements of
$H_{r}(K)$ are equivalence classes in $\cup H_{r}(\mathcal{O})$, where the
union is taken over all the open sets $\mathcal{O}$ containing $K$. Two
elements of $\cup H_{r}(\mathcal{O})$ are said to be equivalent if they
coincide in some open set containing $K$.

For each open $\mathcal{O}\supset K$ let $\rho_{\mathcal{O}}:H_{r}%
(\mathcal{O})\rightarrow H_{r}(K)$ be the map such that $\rho_{\mathcal{O}%
}(u)$ is the equivalence class of $u$. The inductive limit topology is the
finest topology in $H_{r}(K)$ such that all the maps $\rho_{\mathcal{O}}$ are continuous.

If $\mathcal{O}_{1}\subset\mathcal{O}_{2}$ and $u\in H_{r}(\mathcal{O}_{2})$
then $\rho_{\mathcal{O}_{2}}(u)=\rho_{\mathcal{O}_{1}}(\left.  u\right\vert
_{\mathcal{O}_{1}})$. In some cases, we shall write simply $\rho(u)$ instead
of $\rho_{\mathcal{O}}(u)$.

Similar maps corresponding to spaces $Wh_{r}(K),Wh(K)$ and $P_{r}(K)$ will be
denoted by also by $\rho_{\mathcal{O}}$.

Note that all the spaces $H_{r}(K),Wh_{r}(K),Wh(K)$ and $P_{r}(K)$ are locally convex.

We now introduce two more spaces taking in consideration in $P_{r,0}%
(E\setminus K)$ and $P_{r}(K)$ the closed subspaces formed by the elements of
the form $dh_{1}+\delta h_{2}$ and $\rho\left(  dh_{1}+\delta h_{2}\right)  $,
respectively, where $h_{1}$ is an $\left(  r-2\right)  $-form and $h_{2}$ is
an $\left(  r+2\right)  $-form, both of class $C^{2}$ and satisfying
\[
\delta dh_{1}=0,\quad d\delta h_{2}=0
\]
in their respective domains of definition. When $r=1$ we take $h_{1}=0$ and
when $r=n-1$ we take $h_{2}=0$. Thus, the considered spaces are the
corresponding quotient spaces, which will be denoted by $\bar{P}%
_{r,0}(E\setminus K)$ and $\bar{P}_{r}(K)$, respectively. Both these spaces
are provided with the natural quotient topology.

The main goal of this paper is to prove the following assertions:

1) $\left(  H_{r}(K)\right)  ^{\ast}$ is isomorphic to $\bar{P}_{r,0}%
(E\setminus K);$

2) $\left(  \bar{P}_{r}(K)\right)  ^{\ast}$ is isomorphic to $H_{r,0}%
(E\setminus K)$.

As it may be seen through this paper the first assertion plays an important
roll in the formulation of the second one.

Both 1) and 2) are the generalization of the results concerning the dual of
the space of analytic functions on a planar compact set (see \cite{koe}) to
higher dimensions.

\bigskip
\textbf{B.}
It is easy to see that $P_{r}(K)$ is continuously imbedded in $Wh_{r-1}%
(K)+Wh_{r+1}(K)$. Indeed, from the commutativity of the diagram%
\[%
\begin{tabular}
[c]{lll}%
$P_{r}(K)$ & $\overset{i_{1}}{\longrightarrow}$ & $Wh_{r-1}(K)+Wh_{r+1}(K)$\\
$\rho_{\mathcal{O}}\uparrow$ &  & $\rho_{\mathcal{O}}\uparrow$\\
$P_{r}(\mathcal{O})$ & $\overset{i_{2}}{\longrightarrow}$ & $Wh_{r-1}%
(\mathcal{O})+Wh_{r+1}(\mathcal{O})$%
\end{tabular}
\]
where $i_{1},i_{2}$ are the canonical imbeddings and $\mathcal{O}$ is an open
set containing $K$, it follows the continuity of $i_{1}$.

Let $\mathcal{N}$ be the subspace of $Wh_{r-1}(K)+Wh_{r+1}(K)$ whose elements
are of the form $\rho\left(  dh_{1}+\delta h_{2}\right)  $, where $h_{1}%
,h_{2}$ are differential forms such that $dh_{1}\in Wh_{r-1}(\mathcal{O})$ and
$\delta h_{2}\in Wh_{r+1}(\mathcal{O})$ for some open $\mathcal{O}$ containing
$K$. If $r=1$ or $r=n-1$ the space $\mathcal{N}$ just consists of the elements
of the form $\rho\left(  \delta h_{2}\right)  $ or $\rho\left(  dh_{1}\right)
$, respectively. Then, it is clear that $P_{r}(K)\cap\mathcal{N}$ is a
subspace of $P_{r}(K)$, $\bar{P}_{r}(K)=P_{r}(K)\diagup\left(  P_{r}%
(K)\cap\mathcal{N}\right)  $ is continuously imbedded in $\left(
Wh_{r-1}(K)+Wh_{r+1}(K)\right)  \diagup\mathcal{N}$. Note that this latter
space is locally convex.

\bigskip
{\textbf{C.}\label{1.D}} The linear functionals on certain spaces of
differential forms admit to be exteriorly multiplied by differential
forms from those spaces. Let us consider those spaces.

In this subsection the topology of the space does not play any roll. Let
$\chi\left(  \mathcal{O}\right)  $ be linear space of nonhomogeneous
differential forms on a set $\mathcal{O}$ (not necessarily open). Assume that
$\chi\left(  \mathcal{O}\right)  $ satisfies the conditions%
\begin{equation}
\text{if }u\in\chi\left(  \mathcal{O}\right)  \text{ then }\ast u\in
\chi\left(  \mathcal{O}\right)  \text{ and }u\wedge\xi\in\chi\left(
\mathcal{O}\right)  \text{ for each }\xi\in%
{\textstyle\bigwedge^{\ast}}
E. \label{eq02}%
\end{equation}
If $\chi\left(  \mathcal{O}\right)  $ satisfies (\ref{eq02}) it is valid that
$\chi\left(  \mathcal{O}\right)  =%
{\textstyle\bigoplus_{r=0}^{n}}
\chi_{r}\left(  \mathcal{O}\right)  $, where
\[
\chi_{r}\left(  \mathcal{O}\right)  =\left\{  u\in\chi\left(  \mathcal{O}%
\right)  :\text{Im}u\subset%
{\textstyle\bigwedge^{r}}
E\right\}  .
\]
Thus, if $\Lambda$ is a linear functional on $\chi\left(  \mathcal{O}\right)
$ then $\Lambda=\sum_{r=0}^{n}\Lambda_{r}$, where $\Lambda_{r}$ is a linear
functional that vanishes on each subspace $\chi_{s}\left(  \mathcal{O}\right)
$ with $s\neq r$.

Clearly, each linear functional on $\chi_{r}\left(  \mathcal{O}\right)  $ may
be identified with one such $\Lambda_{r}$.

The exterior product $\Lambda\wedge u$ may defined for every linear functional
$\Lambda$ on $\chi\left(  \mathcal{O}\right)  $ and every form $u\in
\chi\left(  \mathcal{O}\right)  $ as a bilinear map with values in $%
{\textstyle\bigwedge^{\ast}}
E$, such that, if $\Lambda=\Lambda_{r}$ and $u\in\chi_{s}\left(
\mathcal{O}\right)  $ then $\Lambda\wedge u$ is the $\left(  r+s\right)
$-covector that satisfies%
\[
\left\langle \Lambda_{r}\wedge u,\xi\right\rangle =\left(  -1\right)
^{nr+n}\Lambda_{r}\left(  \ast\left(  u\wedge\ast\xi\right)  \right)
\]
for every $\xi\in%
{\textstyle\bigwedge^{r+s}}
E$.

When $u\in\chi_{0}\left(  \mathcal{O}\right)  $ the product $\Lambda_{r}\wedge
u$ will be denoted by $\Lambda_{r}u$.

On the other hand, for each linear functional $\Lambda$ on $\chi\left(
\mathcal{O}\right)  $ it is possible to define $\ast\Lambda$ by the equality%
\[
\ast\Lambda:=\sum_{r=0}^{n}\ast\Lambda_{r},
\]
where $\ast\Lambda_{r}$ is the linear functional on $\chi\left(
\mathcal{O}\right)  $ such that $\ast\Lambda_{r}(u)=\left(  -1\right)
^{nr+r}\Lambda_{r}\left(  \ast u\right)  $.

The following properties follow from the definition:

1) $\Lambda_{r}(u)=\ast\left(  \Lambda_{r}\left(  \ast u\right)  \right)  $
for every $u\in\chi\left(  \mathcal{O}\right)  $.

2) $\Lambda\wedge\left(  u\wedge\xi\right)  =\left(  \Lambda\wedge u\right)
\wedge\xi$ for every linear functional $\Lambda$ on $\chi\left(
\mathcal{O}\right)  $, $u\in\chi\left(  \mathcal{O}\right)  $ and $\xi\in%
{\textstyle\bigwedge^{\ast}}
E$.

3) $\ast\left(  \Lambda_{r}\left(  \ast u\right)  \right)  =\left(  -1\right)
^{nr+r}\ast\left(  \ast\Lambda_{r}\wedge u\right)  $ for each $u\in\chi
_{r}\left(  \mathcal{O}\right)  $.

One example of space satisfying (\ref{eq02}) is $Wh(\mathcal{O})$. Thus,
linear functionals on $Wh(\mathcal{O})$ satisfy the equalities 1), 2) and 3).

\bigskip
\textbf{D.} Each linear functional $\Lambda_{r}$ admits a coordinate
expression. Let $e_{1},\cdots,e_{n}$ be an orthonormal basis in $E$
and $e^{1},\cdots,e^{n}$ its dual basis.. Then there exist linear
functionals $\Lambda_{0,\alpha}$ on $Wh_{0}(\mathcal{O})$ (the space
of harmonic functions on $\mathcal{O}$) such
that%
\[
\Lambda_{r}=\sum_{\alpha\in\mathcal{M}_{r}}\Lambda_{0,\alpha}e^{\alpha},
\]
where $\mathcal{M}_{r}$ is the set of all increasing multi-indices $\left(
\alpha_{1},\cdots,\alpha_{r}\right)  $ with $\alpha_{i}\in\left\{
1,\cdots,n\right\}  $. The last equality means that%
\[
\Lambda_{r}(u)=\sum_{\alpha\in\mathcal{M}_{r}}\Lambda_{0,\alpha}(u_{\alpha
})\text{ \ \ for every }u=\sum_{\alpha\in\mathcal{M}_{r}}u_{\alpha}e_{\alpha
}\in\chi_{r}\left(  \mathcal{O}\right)  .
\]

In fact, to obtained the above given coordinate decomposition it is enough to
define $\Lambda_{0,\alpha}(f):=\Lambda_{r}(fe_{\alpha})$ for every $f\in
\chi_{0}\left(  \mathcal{O}\right)  $. Now it is clear that%
\[
\Lambda_{r}\wedge u=\sum_{\alpha\in\mathcal{M}_{r},\beta\in\mathcal{M}_{s}%
}\Lambda_{0,\alpha}(u_{\beta})e^{\alpha}\wedge e^{\beta}\text{ }%
\]
for each $u\in\chi_{s}\left(  \mathcal{O}\right)  $. Note that this formula is
similar to the coordinate expression for the exterior product of differential
forms. In particular, if $f\in\chi_{0}\left(  \mathcal{O}\right)  $ we have%
\[
\Lambda_{r}f:=\Lambda_{r}\wedge f=\sum_{\alpha\in\mathcal{M}_{r}}%
\Lambda_{0,\alpha}(f)e^{\alpha}.
\]

\section{Currents, covectorial charges and Newtonian potential of currents.}

In \cite{pre} it is given a more complete treatment of the facts in this
section. We include here only the most significant details we use in this paper.

\subsection{Exterior product of a current and a form}

Let $\mathcal{D}_{r}$ denote the space of infinitely differentiable
differential forms of degree $r$ with compact support. Let $\mathcal{D}%
_{r}^{\prime}$ denote the space of currents of degree $n-r$ and dimension $r$.
The symbol $\mathcal{E}_{r}^{\prime}$ will be use to represent the subspace of
currents of $\mathcal{D}_{r}^{\prime}$ with compact support.

For each current $T\in\mathcal{D}_{r}^{\prime}$ the symbols $dT,\ast T,\delta
T$ and $\Delta T$ denote currents belonging to $\mathcal{D}_{r-1}^{\prime
},\mathcal{D}_{n-r}^{\prime},\mathcal{D}_{r+1}^{\prime}$ and $\mathcal{D}%
_{r}^{\prime}$, respectively (see \cite{rha} for more details).

Following \cite{rha} we define the exterior product $T\wedge u$ of a current
$T\in\mathcal{E}_{r}^{\prime}$ and an infinitely differentiable $s$-form $u$
defined in a neighborhood of the support of $T$ as a current in $\mathcal{D}%
_{r-s}$ satisfying%
\[
T\wedge u[\omega]=T[u\wedge\omega]
\]
for $\omega\in\mathcal{D}_{r-s}.$

\subsection{Currents determined by a surface}

Now we define two special kinds of currents that will be used in the following
sections. With any smooth oriented hypersurface $\mathcal{S}$ in $E$ (i.e., an
$(n-1)$-dimensional oriented $C^{1}$-submanifold of $E$) \ we associate the
$\left(  n-1\right)  $-dimensional current, denoted again by $\mathcal{S}$,
given by the formula%
\[
\mathcal{S}[\omega]:=\ast\int_{\mathcal{S}}\left(  \bar{N}\wedge\omega\right)
dS
\]
for $\omega\in\mathcal{D}_{n-1}$, where $S$ is the $\left(  n-1\right)
$-dimensional measure on $\mathcal{S}$ and $\bar{N}$ is the field of unitary
normal vectors to $\mathcal{S}$ that determines its orientation.

Moreover, every Lebesgue measurable subset $A$ of $\mathcal{S}$ gives rise to
an $n$ -dimensional current $\chi_{A}$%
\[
A[\omega]:=\int c_{A}\omega=\int_{A}\omega=\ast\int_{A}\omega d\nu,
\]
where $\nu$ is the Lebesgue measure in $E$. More details on the integration of
differential forms with respect to a measure are found in \cite{pre}.

\subsection{Currents determined by a measure}

Let $\mu$ is an $r$-covector valued Borel charge in $E$ (briefly, $\mu\in
M^{r}$). For each $\mu\in M^{r}$ with compact support in an open set
$\mathcal{O}$ in $E$ and a continuous differential form of degree $s$ defined
in $\mathcal{O}$ it is possible to define $\int d\mu\wedge u$ as limit of
Riemann sums $\sum\mu\left(  A_{k}\right)  \wedge u\left(  \xi_{k}\right)  ,$
where $A_{k}$ are Borel sets and $\xi_{k}\in A_{k}$, as $\max\left(
\text{diam }A_{k}\right)  $ tends to zero. In this case, $\int d\mu\wedge u$
exists and is an $\left(  r+s\right)  $-covector. It is obvious that if $u$ is
a continuous $s$-form such that $\int d\mu\wedge u=0$ for every $\mu\in
M^{n-s}$ then $u=0$.

Also, with $\mu\in M^{r}$ we can associate a current, denoted again by $\mu$,
such that%
\begin{equation}
\mu\left[  \omega\right]  :=\ast\int d\mu\wedge\omega,\label{eq03}%
\end{equation}
with $\omega\in\mathcal{D}_{n-r}$.

If $\mu\in M^{r}$ and supp $\mu$ is compact then the current $\mu$ defined by
(\ref{eq03}) belongs to $\mathcal{E}_{n-r}^{\prime}$.

\subsection{Newtonian potential of a current}

The Newtonian potential of a current $T\in\mathcal{E}_{r}^{\prime}$ is defined
by the formula%
\[
U^{T}\left[  \varphi\right]  :=T\left[  U^{\varphi}\right]  ,\quad\varphi
\in\mathcal{D}_{r},
\]
where $U^{\varphi}$ is the Newtonian potential of $\varphi$. From this
definition follows that $U^{T}\in\mathcal{D}_{r}^{\prime}$ (see \cite{pre},
2.1.A and 2.1.B for more details).

Under the given assumptions on $T$ the current $U^{T}$ coincides with a
differential form of class $C^{\infty}$ outside of the support of $T$. Hence,
if $T_{1}\in\mathcal{E}_{r}^{\prime}$, $T_{2}\in\mathcal{E}_{n-r}^{\prime}$
and supp $T_{1}\cap$supp $T_{2}=\varnothing$, both $T_{1}\left[  U^{T_{2}%
}\right]  $ and $T_{2}\left[  U^{T_{1}}\right]  $ have sense and it holds the
reprocity law%
\[
T_{1}\left[  U^{T_{2}}\right]  =\left(  -1\right)  ^{nr+r}T_{2}\left[
U^{T_{1}}\right]  .
\]

\section{Cauchy-Green representation formulas for harmonic forms and
holomorphic pairs}

The Cauchy-Green formulas for harmonic forms and holomorphic pairs may be
considered as generalizations of the classical Cauchy formula for analytic
functions. As in the classic case, these formulas play an important roll in
the description of linear functionals on the space of germs of harmonic forms
on a compact in an Euclidean space.

\subsection{Cauchy-Green representation formulas for harmonic forms}

Let $u$ be an $r$-form harmonic in some open $\mathcal{O}\subset E$ and
$K\subset E$ be a regular compact (i.e., $K$ is a smooth compact manifold with
boundary $\partial K$).

Assuming that $\partial K$ is oriented by means of the exterior normal field
the following formula is true%
\begin{equation}
u=-\frac{1}{c_{c}}\left(  \delta U^{\partial K\wedge u}+\gamma_{r}%
dU^{\ast\left(  \partial K\wedge\ast u\right)  }\right)  \label{eq04}%
\end{equation}
in $\mathring{K}$, where $\gamma_{r}:=\left(  -1\right)  ^{nr+n+1}$ and
$c_{n}$ is the product of the $\left(  n-1\right)  $-dimensional measure of
the sphere $S^{n-1}\subset%
\mathbb{R}
^{n}$ by $n-2$ (see \cite{pre}).

This formula may be rewritten as follows%
\begin{equation}
u(x)=-\frac{1}{c_{n}}\left\{  \delta_{x}\int_{\partial K}\frac{\bar
{N}(y)\wedge u(y)}{\left\Vert x-y\right\Vert ^{n-2}}dS(y)+\gamma_{r}d_{x}%
\ast\int_{\partial K}\frac{\bar{N}(y)\wedge\ast u(y)}{\left\Vert
x-y\right\Vert ^{n-2}}dS(y).\right\}  \label{eq05}%
\end{equation}

\begin{remark}
If $K$ is a regular compact such that $E\smallsetminus K\subset\mathcal{O}$
and $u$ is a harmonic in $\mathcal{O}$ $r$-form that vanishes at infinity,
then a formula similar to (\ref{eq05}), which differs just in the sign, is
valid for $u$ in $E\smallsetminus K$. The proof of this fact is similar to the
classical one for analytic functions of complex variable.
\end{remark}

\subsection{Cauchy-Green representation formulas for holomorphic pairs}

Let $w=w_{r-1}+w_{r+1}$ be a holomorphic pair defined some open $\mathcal{O}%
\subset E$ and $K\subset E$ be a regular compact, then%
\begin{align}
w_{r+1}  &  =-\frac{1}{c_{n}}\left\{  \delta U^{\partial K\wedge w_{r+1}%
}+\gamma_{r+1}dU^{\ast\left(  \partial K\wedge\ast w_{r+1}\right)
}+dU^{\partial K\wedge w_{r-1}}\right\}  ,\label{eq06}\\
w_{r-1}  &  =-\frac{1}{c_{n}}\left\{  \delta U^{\partial K\wedge w_{r-1}%
}+\gamma_{r-1}dU^{\ast\left(  \partial K\wedge\ast w_{r-1}\right)  }%
+\gamma_{r+1}\delta U^{\ast\left(  \partial K\wedge w_{r+1}\right)  }\right\}
, \label{eq07}%
\end{align}
in $\mathring{K}$. The proof of this formula is similar to the proof of
(\ref{eq04}) given in \cite{pre}.

As in the case of harmonic forms, these formulas may be written in terms of
integrals. Moreover, if $E\smallsetminus K\subset\mathcal{O}$ and $w$ is an
holomorphic pair defined in $\mathcal{O}$ vanishing at infinity, then $w$ may
be represented in $E\smallsetminus K$ by a formula that differes from
(\ref{eq06}-\ref{eq07}) just by the sign.

\section{Dual space of $H_{r}(K)$}

In this section we provide a description of the continuous linear functionals
on $H_{r}(K)$.

\begin{theorem}
\label{theo01}Let $K$ be a compact set and $\Lambda:H_{r}(K)\rightarrow%
\mathbb{R}
$. Then, $\Lambda\in\left(  H_{r}(K)\right)  ^{\ast}$ if and only if there
exists a holomorphic pair $w=w_{r-1}+w_{r+1}\in P_{0,r}(E\smallsetminus K)$
such that
\begin{align}
\Lambda\left(  \rho_{\mathcal{O}}(u)\right)  =
&-\frac{1}{c_{n}}\int_{\partial K_{1}}\ast\left(
w_{r+1}\wedge\ast\left(  \bar{N}\wedge u\right)  \right) dS
\nonumber
\\
&+\frac{\left(  -1\right)  ^{r+1}}{c_{n}}\int_{\partial
K_{1}}\ast\left(
w_{r-1}\wedge\left(  \bar{N}\wedge\ast u\right)  \right)  dS \label{eq08}%
\end{align}
for every open neiborhood $\mathcal{O}$ of $K$, any form $u\in H_{r}%
(\mathcal{O})$ and any regular compact $K_{1}$ satisfying $K\subset
\mathring{K}_{1}\subset K_{1}\subset\mathcal{O}$. Moreover, for every
$\Lambda\in\left(  H_{r}(K)\right)  ^{\ast}$ all the pairs $w$ with this
property belongs to the same class in $\bar{P}_{0,r}(E\smallsetminus K)$.
\end{theorem}

\begin{lemma}
\label{lem01}Let $\Omega$ be an open set and $u\in H_{r}(\Omega),w_{r-1}%
+w_{r+1}\in P_{r}(\Omega)$. Then for every regular compact $L\subset\Omega$%
\[
-\int_{\partial K_{1}}\ast\left(  w_{r+1}\wedge\ast\left(  \bar{N}\wedge
u\right)  \right)  dS+\left(  -1\right)  ^{r+1}\int_{\partial K_{1}}%
\ast\left(  w_{r-1}\wedge\left(  \bar{N}\wedge\ast u\right)  \right)  dS=0.
\]

\end{lemma}

\begin{proof}
The proof is based on the Gauss-Ostrogradskii formula (see \cite{pre})%
\[
\int_{L}\ast\left(  du\right)  d\upsilon=\int_{\partial L}\ast\left(  \bar
{N}\wedge u\right)  dS.
\]
The left member of the equality asserted in the lemma may be transformed in
the following way%
\begin{align*}
&  -\int_{\partial K_{1}}\ast\left(  w_{r+1}\wedge\ast\left(  \bar{N}\wedge
u\right)  \right)  dS+\left(  -1\right)  ^{r+1}\int_{\partial K_{1}}%
\ast\left(  w_{r-1}\wedge\left(  \bar{N}\wedge\ast u\right)  \right)  dS\\
&  =-\int_{\partial K_{1}}\ast\left(  \bar{N}\wedge u\wedge\ast w_{r+1}%
\right)  dS-\gamma_{r}\int_{\partial K_{1}}\ast\left(  \bar{N}\wedge\ast
u\wedge w_{r-1}\right)  dS\\
&  =-\int_{L}\ast d\left(  u\wedge\ast w_{r+1}\right)  d\upsilon-\gamma
_{r}\int_{L}\ast d\left(  \ast u\wedge w_{r-1}\right)  d\upsilon.
\end{align*}
Since $u\in H_{r}(\Omega)$ and $w_{r-1}+w_{r+1}\in P_{r}(\Omega)$ we have that
$d\left(  u\wedge\ast w_{r+1}\right)  +\gamma_{r}d\left(  \ast u\wedge
w_{r-1}\right)  =0$. So, the lemma is proved.
\end{proof}

\begin{proof}
[Proof of the theorem]\textbf{First part}. We shall prove here that formula
(\ref{eq08}) defines a continuous linear functional on $H_{r}(K)$. At first,
let us see that $\Lambda$ is well defined. The right hand member of
(\ref{eq08}) depends only on $\rho_{\mathcal{O}}(u)$. The independence on $u$
and $K_{1}$ is deduced from the Lemma \ref{lem01}.

Let us note that the right hand term of (\ref{eq08}) depends just on the class
whose representative element is $w_{r-1}+w_{r+1}$ and not on the particular
representative element. It is due to the fact that, if we consider
$w_{r-1}=dh_{1}$ and $w_{r+1}=\delta h_{2}$ (see the definition of $\bar
{P}_{0,r}(E\smallsetminus K)$) it follows from the Stokes formula that the
integral in (\ref{eq08}) vanishes.

Obviously, $\Lambda$ is linear. In order to prove that $\Lambda$ is continuous
it is enough to show that $\Lambda\circ\rho_{\mathcal{O}}$ is continuous for
every open neighborhood $\mathcal{O}$ of $K$.

Let $K_{1}$ be a regular compact such that $K\subset\mathring{K}_{1}\subset
K_{1}\subset\mathcal{O}$, then%
\begin{equation}
\left\vert \Lambda\left(  \rho_{\mathcal{O}}(u)\right)  \right\vert
=\left\vert \Lambda\left(  \rho_{\mathring{K}_{1}}(u)\right)  \right\vert \leq
C\max_{x\in\partial K_{1}}\left\vert u(x)\right\vert , \label{eq09}%
\end{equation}
where $C$ is a positive constant independent of $u$. To obtain inequality
(\ref{eq09}) we have used the fact that%
\[
\left\vert \xi\wedge\eta\right\vert \leq\binom{r+s}{r}\left\vert
\xi\right\vert \left\vert \eta\right\vert
\]
for $\xi\in%
{\textstyle\bigwedge^{r}}
E$ and $\eta\in%
{\textstyle\bigwedge^{s}}
E$ (see, e.g., \cite{fed}).

From (\ref{eq09}) it follows the continuity of $\Lambda\circ\rho_{\mathcal{O}%
}$ and thus the continuity of $\Lambda$.

\textbf{Second part}. Let us show now that every $\Lambda\in\left(
H_{r}(K)\right)  ^{\ast}$ may be represented by (\ref{eq08}) for some
$w_{r-1}+w_{r+1}\in P_{0,r}(E\smallsetminus K)$. Let $\mathcal{O}$ be an open
neiborhood of $K$, $u\in H(\mathcal{O})$ and $K_{1}$ a regular compact
satisfying $K\subset\mathring{K}_{1}\subset K_{1}\subset\mathcal{O}$. Using
the Cauchy-Green formula (\ref{eq05}) we have that $u(x)=u_{1}(x)+u_{2}(x)$
where%
\begin{align*}
u_{1}(x)  &  =-\frac{1}{c_{n}}\delta_{x}\int_{\partial K_{1}}\frac{\bar
{N}(y)\wedge u(y)}{\left\Vert x-y\right\Vert ^{n-2}}dS(y)\\
&  =-\frac{\gamma_{r+1}}{c_{n}}\int_{\partial K_{1}}\ast\left(  d_{x}\frac
{1}{\left\Vert x-y\right\Vert ^{n-2}}\wedge\ast\left(  \bar{N}(y)\wedge
u(y)\right)  \right)  dS(y),
\end{align*}%
\begin{align*}
u_{2}(x)  &  =-\frac{\gamma_{r}}{c_{n}}d_{x}\ast\int_{\partial K}\frac{\bar
{N}(y)\wedge\ast u(y)}{\left\Vert x-y\right\Vert ^{n-2}}dS(y)\\
&  =-\frac{\gamma_{r}}{c_{n}}\int_{\partial K}\left(  d_{x}\frac{1}{\left\Vert
x-y\right\Vert ^{n-2}}\wedge\ast\left(  \bar{N}(y)\wedge\ast u(y)\right)
\right)  dS(y).
\end{align*}
Moreover, since $H_{r}(K)$ is continuously imbedded in $Wh_{r}(K)$, for each
$\Lambda\in\left(  H_{r}(K)\right)  ^{\ast}$ there exists a functional
$\Lambda^{\prime}\in\left(  Wh_{r}(K)\right)  ^{\ast}$ that extends $\Lambda$.

Due to the reasons given above and to the continuity of $\Lambda^{\prime}%
\circ\rho_{\mathcal{O}}$ on $Wh_{r}(\mathring{K}_{1})$ it is valid that
$\Lambda^{\prime}\circ\rho_{\mathcal{O}}$ commutes with the integral sign.
Thus, using the properties 1) and 2) from subsection \ref{1.D} we obtain that%
\begin{align}
\Lambda^{\prime}\circ\rho_{\mathcal{O}}(u_{1})  &  =-\frac{\gamma_{r+1}}%
{c_{n}}\int_{\partial K_{1}}\ast\left\{  \Lambda^{\prime}\circ\rho
_{\mathcal{O}}\wedge\ast\ast\left(  d_{x}\frac{1}{\left\Vert x-y\right\Vert
^{n-2}}\wedge\ast\left(  \bar{N}(y)\wedge u(y)\right)  \right)  \right\}
dS(y)\label{eq10}\\
&  =\frac{(-1)^{r}}{c_{n}}\int_{\partial K_{1}}\ast\left\{  \left(
\Lambda^{\prime}\circ\rho_{\mathcal{O}}\wedge d_{x}\frac{1}{\left\Vert
x-y\right\Vert ^{n-2}}\right)  \wedge\ast\left(  \bar{N}(y)\wedge u(y)\right)
\right\}  dS(y).\nonumber
\end{align}

Let us remark that
\[
d_{x}\frac{1}{\left\Vert x-y\right\Vert ^{n-2}}=-d_{y}\frac{1}{\left\Vert
x-y\right\Vert ^{n-2}}%
\]
and that, because of the continuity of $\Lambda^{\prime}\circ\rho
_{\mathcal{O}}$, the form $\Lambda^{\prime}\circ\rho_{\mathcal{O}}\frac
{1}{\left\Vert x-y\right\Vert ^{n-2}}$ (depending on $y$) is of class $C^{1}$.
Thus it can be proved that%
\begin{equation}
\Lambda^{\prime}\circ\rho_{\mathcal{O}}\wedge d_{x}\frac{1}{\left\Vert
x-y\right\Vert ^{n-2}}=-\Lambda^{\prime}\circ\rho_{\mathcal{O}}\wedge
d_{y}\frac{1}{\left\Vert x-y\right\Vert ^{n-2}}=(-1)^{r+1}d_{y}\left(
\Lambda^{\prime}\circ\rho_{\mathcal{O}}\frac{1}{\left\Vert x-y\right\Vert
^{n-2}}\right)  . \label{eq11}%
\end{equation}
Indeed, the last equality if verified as follows. Let $\varphi\in
\mathcal{D}_{r+1}$ with supp $\varphi\subset E\smallsetminus K$ then,%
\begin{align*}
\int\left\langle d_{y}\left(
\Lambda^{\prime}\circ\rho_{\mathcal{O}}\frac {1}{\left\Vert
x-y\right\Vert ^{n-2}}\right)  ,\varphi(y)\right\rangle
d\upsilon(y)\\=\int\left\langle
\Lambda^{\prime}\circ\rho_{\mathcal{O}}\frac {1}{\left\Vert
x-y\right\Vert ^{n-2}},\delta\varphi(y)\right\rangle d\upsilon(y)
\\
=(-1)^{r+1}\int\left(  \Lambda^{\prime}\circ\rho_{\mathcal{O}}\wedge
d_{y}\left(  \frac{\ast\varphi(y)}{\left\Vert x-y\right\Vert
^{n-2}}\right) \right)  d\upsilon(y)\\+(-1)^{r}\int\left(
\Lambda^{\prime}\circ\rho _{\mathcal{O}}\wedge
d_{y}\frac{1}{\left\Vert x-y\right\Vert ^{n-2}}\wedge
\ast\varphi(y)\right)  d\upsilon(y)=
\\
=(-1)^{r+1}\ast\left(
\Lambda^{\prime}\circ\rho_{\mathcal{O}}\wedge\int d_{y}\left(
\frac{\ast\varphi(y)}{\left\Vert x-y\right\Vert ^{n-2}}\right)
d\upsilon(y)\right) \\ +(-1)^{r}\int\left\langle
\Lambda^{\prime}\circ
\rho_{\mathcal{O}}\wedge d_{y}\frac{1}{\left\Vert x-y\right\Vert ^{n-2}%
},\varphi(y)\right\rangle d\upsilon(y)=
\\
=(-1)^{r}\int\left\langle \Lambda^{\prime}\circ\rho_{\mathcal{O}}\wedge
d_{y}\frac{1}{\left\Vert x-y\right\Vert ^{n-2}},\varphi(y)\right\rangle
d\upsilon(y).
\end{align*}
In the last two equalities we have used the facts that $\Lambda^{\prime}%
\circ\rho_{\mathcal{O}}$ is continuous and that%
\[
\int d_{y}\left(  \frac{\ast\varphi(y)}{\left\Vert x-y\right\Vert ^{n-2}%
}\right)  d\upsilon(y)=0
\]
for every $x\notin$ supp $\varphi$.

Thus, relations (\ref{eq10}) and (\ref{eq11}) imply that%
\[
\Lambda^{\prime}\circ\rho_{\mathcal{O}}(u_{1})=-\frac{1}{c_{n}}\int_{\partial
K_{1}}\ast\left(  w_{r+1}\wedge\ast\left(  \bar{N}\wedge u\right)  \right)
dS,
\]
where%
\[
w_{r+1}=d_{y}\left(  \Lambda^{\prime}\circ\rho_{\mathcal{O}}\frac
{1}{\left\Vert x-y\right\Vert ^{n-2}}\right)  .
\]

Similarly,
\begin{align*}
\Lambda\circ\rho_{\mathcal{O}}(u_{2})  &  =\Lambda^{\prime}\circ
\rho_{\mathcal{O}}(u_{2})\\&=-\frac{\gamma_{r}}{c_{n}}\int_{\partial K_{1}}%
\ast\left\{  \Lambda^{\prime}\circ\rho_{\mathcal{O}}\wedge\ast\left(
d_{x}\frac{1}{\left\Vert x-y\right\Vert ^{n-2}}\wedge\ast\left(  \bar
{N}(y)\wedge\ast u(y)\right)  \right)  \right\}  dS(y)\\
&  =-\frac{\gamma_{r}}{c_{n}}\int_{\partial K_{1}}\ast\left\{  \ast
\Lambda^{\prime}\circ\rho_{\mathcal{O}}\wedge\ast\ast\left(  d_{x}\frac
{1}{\left\Vert x-y\right\Vert ^{n-2}}\wedge\ast\left(  \bar{N}(y)\wedge\ast
u(y)\right)  \right)  \right\}  dS(y)\\
&  =-\frac{1}{c_{n}}\int_{\partial K_{1}}\ast\left\{  d_{y}\ast\left(
\Lambda^{\prime}\circ\rho_{\mathcal{O}}\frac{1}{\left\Vert x-y\right\Vert
^{n-2}}\right)  \wedge\ast\left(  \bar{N}(y)\wedge\ast u(y)\right)  \right\}
dS(y)\\
&  =-\frac{1}{c_{n}}\int_{\partial K_{1}}\ast\left\{  \ast d_{y}\ast\left(
\Lambda^{\prime}\circ\rho_{\mathcal{O}}\frac{1}{\left\Vert x-y\right\Vert
^{n-2}}\right)  \wedge\ast\ast\left(  \bar{N}(y)\wedge\ast u(y)\right)
\right\}  dS(y)\\
&  =\frac{\left(  -1\right)  ^{r+1}}{c_{n}}\int_{\partial K_{1}}\ast\left(
w_{r-1}(y)\wedge\left(  \bar{N}(y)\wedge\ast u(y)\right)  \right)  dS(y),
\end{align*}
where%
\[
w_{r-1}=\delta_{y}\left(  \Lambda^{\prime}\circ\rho_{\mathcal{O}}\frac
{1}{\left\Vert x-y\right\Vert ^{n-2}}\right)  .
\]
Here we have use the property 1) of subsection \ref{1.D}, too.

Thus, the equality (\ref{eq08}) is proved. Now it is easy to see that
$w=w_{r-1}+w_{r+1}\in P_{0,r}(E\smallsetminus K)$. In fact,
\[
dw_{r-1}+\delta w_{r+1}=\Delta\left(  \Lambda^{\prime}\circ\rho_{\mathcal{O}%
}\frac{1}{\left\Vert x-y\right\Vert ^{n-2}}\right)  =0
\]
in $E\smallsetminus K$. To show this last equality it is enough to express the
form $\Lambda^{\prime}\circ\rho_{\mathcal{O}}\frac{1}{\left\Vert
x-y\right\Vert ^{n-2}}$ in coordinates in an orthonormal basis of $E$. The
coefficients of the forms obtained in this way are harmonic functions in
$E\smallsetminus K$.

It just remains to prove that $w_{r-1}(y)\rightarrow0$ and $w_{r+1}%
(y)\rightarrow0$ as $y\rightarrow\infty$ in the norm of $%
{\textstyle\bigwedge^{r}}
E$. Clearly, this is equivalent to $\left\langle w_{r-1}(y),\xi\right\rangle
\rightarrow0$ and $\left\langle w_{r+1}(y),\xi\right\rangle \rightarrow0$ as
$y\rightarrow\infty$ for every $\xi\in%
{\textstyle\bigwedge^{r}}
E$. Note that%
\begin{align*}
\left\langle w_{r+1}(y),\xi\right\rangle =\left(  -1\right)
^{r}\left\langle \Lambda^{\prime}\circ\rho_{\mathcal{O}}\wedge
d_{y}\frac{1}{\left\Vert x-y\right\Vert ^{n-2}},\xi\right\rangle
\\=\pm\Lambda^{\prime}\circ \rho_{\mathcal{O}}\left(  \ast\left(
d_{x}\frac{1}{\left\Vert x-y\right\Vert ^{n-2}}\wedge\ast\xi\right)
\right)  .
\end{align*}
The continuity of $\Lambda^{\prime}\circ\rho_{\mathcal{O}}$ implies that the
last term in the last equality tends to zero as $y\rightarrow\infty$. The
proof for $w_{r-1}$ is completely analogous.

Let us prove now that $\Lambda$ depends only on the class of $P_{0,r}%
(E\smallsetminus K)$ with representative element $w=w_{r-1}+w_{r+1}$ and not
on this particular element. It is enough to show that the periods of both
$\ast w_{r+1}$ and $w_{r-1}$ are uniquely determined by the functional
$\Lambda$. More precisely, we prove the equalities%
\begin{equation}
\int_{\lambda_{1}}\ast w_{r+1}=\Lambda\left(  \rho_{\mathcal{O}}\left(  \delta
U^{\lambda_{1}}\right)  \right)  ,\quad\int_{\lambda_{2}}w_{r-1}%
=\Lambda\left(  \rho_{\mathcal{O}}\left(  d\ast U^{\lambda_{2}}\right)
\right)  ,\label{eq12}%
\end{equation}
for every $(n-r-1)$-dimensional cycle $\lambda_{1}$ and every $(r-1)$%
-dimensional cycle $\lambda_{2}$ contained all in $E\smallsetminus K$. Under
assumed conditions, both $\delta U^{\lambda_{1}}$ and $d\ast U^{\lambda_{2}}$
are harmonic forms in some neighborhood of $K$ (see Pr. 3 in \cite{pre}).
Since these forms are locally integrable they generate currents belonging to
$\mathcal{E}_{r}^{\prime}$.

Let us prove the first equality (\ref{eq12}). By replacing $u=\delta
U^{\lambda_{1}}$ in formula (\ref{eq08}) we obtain%
\begin{align}
\Lambda\left(  \rho_{\mathcal{O}}\left(  \delta
U^{\lambda_{1}}\right) \right)  =-\frac{1}{c_{n}}\int_{\partial
K_{1}}\ast\left(  w_{r+1}\wedge \ast\left(  \bar{N}\wedge\delta
U^{\lambda_{1}}\right)  \right) dS \nonumber\\+\frac{\left(
-1\right) ^{r+1}}{c_{n}}\int_{\partial K_{1}}\ast\left(
w_{r-1}\wedge\left( \bar{N}\wedge\ast\delta U^{\lambda_{1}}\right)
\right)
dS, \label{eq13}%
\end{align}
where $K_{1}$ is a regular compact such that $K\subset\mathring{K}_{1}$ and
$K_{1}\cap$ supp $\lambda_{1}=\varnothing$.

Let us denote by $I$ the first integral in (\ref{eq13}) and by $II$ the second
one. We have%
\begin{align*}
I &  =\frac{\left(  -1\right)  ^{nr+n+r}}{c_{n}}\int_{\partial K_{1}}%
\ast\left(  \ast w_{r+1}\wedge\bar{N}\wedge\delta U^{\lambda_{1}}\right)  dS\\
&  =\frac{\left(  -1\right)  ^{nr+1}}{c_{n}}\int_{\partial K_{1}}\ast\left(
\bar{N}\wedge\ast w_{r+1}\wedge\delta U^{\lambda_{1}}\right)  dS\\
&  =\frac{\left(  -1\right)  ^{nr+1}}{c_{n}}\left(  \partial K_{1}\wedge\ast
w_{r+1}\right)  \left[  \delta U^{\lambda_{1}}\right]  .
\end{align*}
Here $\partial K_{1}\wedge\ast w_{r+1}$ is understood a current belonging to
$\mathcal{E}_{r}^{\prime}$ and thus the above last term has sense. Now, using
the properties given in 2.A and 2.D we have%
\begin{align*}
I &  =\frac{\left(  -1\right)  ^{nr+n+r+1}}{c_{n}}\delta\left(  \partial
K_{1}\wedge\ast w_{r+1}\right)  \left[  U^{\lambda_{1}}\right]  =\frac
{1}{c_{n}}\lambda_{1}\left[  U^{\delta\left(  \partial K_{1}\wedge\ast
w_{r+1}\right)  }\right]  \\
&  =\frac{1}{c_{n}}\lambda_{1}\left[  \delta U^{\partial K_{1}\wedge\ast
w_{r+1}}\right]  =\frac{\gamma_{r+1}}{c_{n}}\lambda_{1}\left[  \ast
dU^{\ast\left(  \partial K_{1}\wedge\ast w_{r+1}\right)  }\right]  .
\end{align*}

By means of similar transformations we may obtain%
\[
II=\frac{1}{c_{n}}\lambda_{1}\left[  \ast dU^{\partial K_{1}\wedge\ast
w_{r-1}}\right]  .
\]
Hence, the first of the equalities (\ref{eq12}) is proved. The proof
of the second equality is similar.

Let us remark that the first equality (\ref{eq12}) implies that $w_{r+1}$ is
uniquely determined up to one term of the form $\delta h_{2}$ where $h_{2}$ is
an $(r+2)$-form such that $d\delta h_{2}=0$. When $r=n-1$ the above condition
simply says that $w_{r+1}$ is uniquely determined. On the other hand, the
second equality (\ref{eq12}) implies that $w_{r-1}$ is determined up to one
term of the form $dh_{1}$, where $h_{1}$ is an $(r-2)$-form such that $\delta
dh_{1}=0$. If $r=1$, simply $w_{r-1}$ is uniquely determined.

Therefore, the proof of the theorem is complete.
\end{proof}

\section{Construction of a space whose dual is the space of harmonic
differential forms}

The results of the previous section allow us to conjecture about which is the
dual space of $H_{r,0}(E\smallsetminus K)$. It seems natural that this space
be $\bar{P}_{r}(K)$. In this section we verify that this fact is indeed true.
More precisely, let $\bar{P}_{r}(K)$ be endowed with the topology induced by
that of $P_{r}(K)$. Then

\begin{theorem}
\label{theo02}Let $K$ be a compact set and $\Lambda:\bar{P}_{r}(K)\rightarrow%
\mathbb{R}
$. Then, $\Lambda\in\left(  \bar{P}_{r}(K)\right)  ^{\ast}$ if and only if
there exists a form $u\in H_{r,0}(E\smallsetminus K)$ such that
\begin{align}
\Lambda\left(  \overline{\rho_{\mathcal{O}}(w)}\right) =\frac{\left(
-1\right)  ^{nr+r+1}}{c_{n}}\int_{\partial K_{1}}\ast\left(
w_{r+1}\wedge \ast\left(  \bar{N}\wedge u\right) \right)
dS\\+\frac{\left(  -1\right) ^{n+1}}{c_{n}}\int_{\partial
K_{1}}\ast\left(  w_{r-1}\wedge\left(  \bar
{N}\wedge\ast u\right)  \right)  dS \label{eq14}%
\end{align}
for every open neiborhood $\mathcal{O}$ of $K$, any par $w=w_{r-1}+w_{r+1}\in
P_{r}(\mathcal{O})$ and any regular compact $K_{1}$ satisfying $K\subset
\mathring{K}_{1}\subset K_{1}\subset\mathcal{O}$.
\end{theorem}

\begin{proof}
\textbf{First part.} In this part it is verified that formula (\ref{eq14})
defines a linear and continuous functional on $\bar{P}_{r}(K)$.

Let us see first that $\Lambda$ given by (\ref{eq14}) is well defined. In view
of Lemma \ref{lem01}, the expression (\ref{eq14}) does not depend neither on
the particular choice of the germ $\rho_{\mathcal{O}}\left(  w\right)  $ nor
on the compact $K_{1}$. To understand that $\Lambda$ does not depend on the
particular germ $\rho_{\mathcal{O}}\left(  w\right)  $ and just on the
equivalence class that it defines it is enough to see that the right hand term
in (\ref{eq14}) vanishes when $\rho_{\mathcal{O}}\left(  w\right)
=\rho_{\mathcal{O}}\left(  dh_{1}+\delta h_{2}\right)  $, where $h_{1}$ is an
$(r-2)$-form and $h_{2}$ an $(r+2)$-form, both of class $C^{2}$ in a
neighborhood of $K$ and satisfying the equations $d\delta h_{1}=0$, $\delta
dh_{2}=0$ in their domains of definition. In fact, it is true due to the fact
that both integrals in (\ref{eq14}) may be written in the form%
\[
\int_{\partial K_{1}}\left(  \bar{N}\wedge d\psi\right)  dS,
\]
where $\psi$ is an $(r-2)$-form, smooth in a neighborhood of $\partial K_{1}$.
Now it suffices to observe that, in view of the Stokes formula, the latter
integral vanishes.

The linearity of $\Lambda$ is immediate. The continuity of $\Lambda$ is
equivalent to the continuity of all the functionals $\Lambda\circ
\rho_{\mathcal{O}}$ for every open $\mathcal{O}$ containing $K$. We have that%
\[
\left\vert \Lambda\left(  \rho_{\mathcal{O}}(w)\right)  \right\vert \leq
C_{1}\max_{x\in\partial K_{1}}\left\vert w_{r-1}(x)\right\vert +C_{2}%
\max_{x\in\partial K_{1}}\left\vert w_{r+1}(x)\right\vert ,
\]
for every regular compact $K_{1}$ satisfying $K\subset\mathring{K}_{1}\subset
K_{1}\subset\mathcal{O}$.

Second part. Let us now prove that for each linear functional $\Lambda
\in\left(  \bar{P}_{r}(K)\right)  ^{\ast}$ there exists a form $u\in
H_{r,0}(E\smallsetminus K)$ that represents $\Lambda$ according to formula
(\ref{eq14}).

Let $\Lambda\in\left(  \bar{P}_{r}(K)\right)  ^{\ast}$. Since
$\bar{P}_{r}(K)$ is continuously imbedded in $$\left(
Wh_{r-1}(K)+Wh_{r+1}(K)\right) \diagup\mathcal{N}$$ and the latter
space is locally convex, there exists $$\Lambda^{\prime}\in\left(
\left(  Wh_{r-1}(K)+Wh_{r+1}(K)\right) \diagup\mathcal{N}\right)
^{\ast}$$ that extends $\Lambda$ and then,
$$\Lambda^{\prime\prime}\in\left(  \left(
Wh_{r-1}(K)+Wh_{r+1}(K)\right)
\right)  ^{\ast}$$ defined by the equality $$\Lambda^{\prime\prime}%
(\rho_{\mathcal{O}}(w))=\Lambda^{\prime}(\rho_{\mathcal{O}}(w))$$
extends $\Lambda$ to $$\left(  Wh_{r-1}(K)+Wh_{r+1}(K)\right).  $$

By definition
$$\Lambda^{\prime\prime}(\rho_{\mathcal{O}}(dh_{1}+\delta
h_{2}))=0$$
for all forms $h_{1},h_{2}$, where $h_{1}$is an $(r-2)$-form and
$h_{2}$ an $(r+2)$-form, both defined in $\mathcal{O}$ and
satisfying $\delta dh_{1}=0,d\delta h_{2}=0$ in $\mathcal{O}$. When
$r=1$ or $r=n-1$, we take $h_{1}=0$ or $h_{2}=0$, respectively.

Let $w$ is represented by formulas (\ref{eq06})-(\ref{eq07}). Thus, in view of
the equality%
\[
\gamma_{r-1}dU^{\ast\left(  \partial K_{1}\wedge\ast w_{r-1}\right)  }+\delta
U^{\partial K_{1}\wedge w_{r+1}}=0,
\]
we have%
\begin{align}
c_{n}\Lambda\left(  \overline{\rho_{\mathcal{O}}(w)}\right)
=-\Lambda ^{\prime\prime}\circ\rho_{\mathcal{O}}\left( \delta
U^{\partial K_{1}\wedge w_{r-1}}+\right.\label{eq15}\\+\left.
\gamma_{r+1}dU^{\ast\left(
\partial K_{1}\wedge\ast w_{r+1}\right) }+\gamma_{r+1}\delta
U^{\ast\left(  \partial K_{1}\wedge w_{r+1}\right) }+dU^{\partial
K_{1}\wedge w_{r-1}}\right)  . \nonumber%
\end{align}

Since $$\Lambda^{\prime\prime}\in\left(  \left(  Wh_{r-1}(K)+Wh_{r+1}%
(K)\right)  \right)  ^{\ast},$$ then there exist
$$\Lambda_{r-1}\in\left( Wh_{r-1}(K)\right)  ^{\ast}, \qquad\Lambda_{r+1}\in\left(  Wh_{r+1}(K)\right)
^{\ast}$$ such that%
\begin{equation}
\Lambda^{\prime\prime}(\rho_{\mathcal{O}}(w_{r-1}+w_{r+1}))=\Lambda_{r-1}%
(\rho_{\mathcal{O}}(w_{r-1}))+\Lambda_{r+1}(\rho_{\mathcal{O}}(w_{r+1})).
\label{eq16}%
\end{equation}

Obviously,
\begin{equation}
\Lambda_{r-1}(\rho_{\mathcal{O}}(dh_{1}))=0,\quad\Lambda_{r+1}(\rho
_{\mathcal{O}}(\delta h_{2})) \label{eq17}%
\end{equation}
if $h_{1}$is an $(r-2)$-form and $h_{2}$ an $(r+2)$-form, both defined in
$\mathcal{O}$ and satisfying $\delta dh_{1}=0,d\delta h_{2}=0$ in
$\mathcal{O}$.

From (\ref{eq15})-(\ref{eq16}) we have that%
\begin{align*}
c_{n}\Lambda\left(  \overline{\rho_{\mathcal{O}}(w)}\right)  =-\Lambda
_{r-1}\circ\rho_{\mathcal{O}}\left(  \delta U^{\partial K_{1}\wedge w_{r-1}%
}+\gamma_{r+1}dU^{\ast\left(  \partial K_{1}\wedge\ast
w_{r+1}\right) }\right) \\
-\Lambda_{r+1}\circ\rho_{\mathcal{O}}\left( \gamma_{r+1}\delta
U^{\ast\left(  \partial K_{1}\wedge w_{r+1}\right)  }+dU^{\partial
K_{1}\wedge w_{r-1}}\right)  .
\end{align*}
Using transformations similar to those done in the proof of
Theorem
\ref{theo01} we obtain that%
\begin{equation}
c_{n}\Lambda\left(  \overline{\rho_{\mathcal{O}}(w)}\right)  =\left(
-1\right)  ^{n+1}\ast\int_{\partial K_{1}}\left(  u\wedge\left(  \bar{N}%
\wedge\ast w_{r+1}+\left(  -1\right)  ^{r+1}\ast\left(  \bar{N}\wedge
w_{r-1}\right)  \right)  \right)  dS, \label{eq18}%
\end{equation}
where%
\[
u(y)=d_{y}\left(  \Lambda_{r-1}\circ\rho_{\mathcal{O}}\frac{1}{\left\Vert
x-y\right\Vert ^{n-2}}\right)  +\delta_{y}\left(  \Lambda_{r+1}\circ
\rho_{\mathcal{O}}\frac{1}{\left\Vert x-y\right\Vert ^{n-2}}\right)
\]
(here $\Lambda_{r-1}\circ\rho_{\mathcal{O}}$ and $\Lambda_{r+1}\circ
\rho_{\mathcal{O}}$ act on form depending on $x$).

Thus, (\ref{eq14}) follows from (\ref{eq18}). Let us prove that in fact $u\in
H_{r,0}(E\smallsetminus K)$. As in the proof of Theorem \ref{theo01} we have
that%
\[
\Delta_{y}\left(  \Lambda_{r-1}\circ\rho_{\mathcal{O}}\frac{1}{\left\Vert
x-y\right\Vert ^{n-2}}\right)  =0,\quad\Delta_{y}\left(  \Lambda_{r+1}%
\circ\rho_{\mathcal{O}}\frac{1}{\left\Vert x-y\right\Vert ^{n-2}}\right)  =0
\]
for $y\in E\smallsetminus K$. Thus, for $y\in E\smallsetminus K$
\[
du(x)=d\delta\left(  \Lambda_{r+1}\circ\rho_{\mathcal{O}}\frac{1}{\left\Vert
x-y\right\Vert ^{n-2}}\right)  =\delta d\left(  \Lambda_{r+1}\circ
\rho_{\mathcal{O}}\frac{1}{\left\Vert x-y\right\Vert ^{n-2}}\right)
\]
and%
\[
\delta u(x)=\delta d\left(  \Lambda_{r-1}\circ\rho_{\mathcal{O}}\frac
{1}{\left\Vert x-y\right\Vert ^{n-2}}\right)  =d\delta\left(  \Lambda
_{r-1}\circ\rho_{\mathcal{O}}\frac{1}{\left\Vert x-y\right\Vert ^{n-2}%
}\right)  .
\]
Consequently, it just remains to prove that, for $y\in E\smallsetminus K$
\begin{equation}
\delta\left(  \Lambda_{r-1}\circ\rho_{\mathcal{O}}\frac{1}{\left\Vert
x-y\right\Vert ^{n-2}}\right)  =0,\quad d\left(  \Lambda_{r-1}\circ
\rho_{\mathcal{O}}\frac{1}{\left\Vert x-y\right\Vert ^{n-2}}\right)  =0.
\label{eq19}%
\end{equation}

Let us prove the first equality, the second one being similar. Let $\xi$ be an
$(r-2)$-covector then, from (\ref{eq17}) we have%
\begin{align*}
\left\langle \delta\left(  \Lambda_{r-1}\circ\rho_{\mathcal{O}}\frac
{1}{\left\Vert x-y\right\Vert ^{n-2}}\right)  ,\xi\right\rangle  &  =\pm
\ast\left(  \ast\Lambda_{r-1}\circ\rho_{\mathcal{O}}\wedge d_{y}\frac
{1}{\left\Vert x-y\right\Vert ^{n-2}}\wedge\ast\xi\right) \\
&  =\pm\ast\left(  \ast\Lambda_{r-1}\circ\rho_{\mathcal{O}}\wedge d_{y}\left(
\frac{\ast\xi}{\left\Vert x-y\right\Vert ^{n-2}}\right)  \right) \\
&  =\pm\Lambda_{r-1}\circ\rho_{\mathcal{O}}\left(  d_{x}\frac{\ast\xi
}{\left\Vert x-y\right\Vert ^{n-2}}\right)  =0
\end{align*}
for $y\in E\smallsetminus K$ \ Therefore, the first equality (\ref{eq19}) is
proved. From the definition of $u$ and the continuity of $\Lambda_{r-1}$ and
$\Lambda_{r+1}$ we have that $u(y)\rightarrow0$ as $y\rightarrow\infty$.

Now it just remains to verify that $u$ is unique. Let $\mu$ be an
$r$-covectorial Borelian charge with support in $E\smallsetminus K$. Then, if
we define $w_{r-1}^{\mu}=\delta U^{\mu}$, $w_{r+1}^{\mu}=dU^{\mu}$ we have
$w^{\mu}=w_{r-1}^{\mu}+w\mu_{r+1}\in P_{r}(E\smallsetminus$supp $\mu)$.
Consequently, substituting in (\ref{eq14}) we obtain%
\begin{align*}
\Lambda\left(  \overline{\rho_{\mathcal{O}}(w^{\mu})}\right)
=\frac{\left( -1\right)  ^{n+r+1}}{c_{n}}\int_{\partial
K_{1}}\ast\left(  dU^{\mu}\wedge \ast\left(  \bar{N}\wedge u\right)
\right)  dS\\+\frac{\left(  -1\right) ^{n+1}}{c_{n}}\int_{\partial
K_{1}}\ast\left(  \delta U^{\mu}\wedge\left( \bar{N}\wedge\ast
u\right)  \right)  dS.
\end{align*}
Proceeding as in the proof of equality (\ref{eq12}) we may proof that
\[
\Lambda\left(  \overline{\rho_{\mathcal{O}}(w^{\mu})}\right)  =\left(
-1\right)  ^{n+r+1}\mu\left[  \ast\left(  \delta U^{\partial K_{1}\wedge
u}+\gamma_{r}dU^{\ast\left(  \partial K_{1}\wedge\ast u\right)  }\right)
\right]  =\left(  -1\right)  ^{n+r+1}\mu\left[  \ast u\right]  .
\]
Therefore, from 2.C if $\Lambda=0$ then $u=0$. The proof of the theorem is complete.
\end{proof}

\subsection{A representation formula for harmonic forms vanishing at the infinity}

As an example of application of the previous results we obtain the
representation formula:

\begin{corollary}
\label{cor01}If $u\in H_{r,0}(E\smallsetminus K)$ then there exist $u_{1}\in
Wh_{r-1}(E\smallsetminus K)$ and $u_{2}\in Wh_{r+1}(E\smallsetminus K)$
vanishing at infinity such that $\delta u_{1}=0$, $du_{2}=0$ and
\begin{equation}
u=du_{1}+\delta u_{2}. \label{eq20}%
\end{equation}

\end{corollary}

\begin{proof}
For every $u\in H_{r,0}(E\smallsetminus K)$ there exists a functional
$\Lambda\in\left(  \bar{P}_{r}(K)\right)  ^{\ast}$ according to (\ref{eq14}).
The second part of the proof of Theorem \ref{theo02} ensures the existence of
a unique form in $H_{r,0}(E\smallsetminus K)$ corresponding to $\Lambda$,
which in addition admits a representation like (\ref{eq20}). Then, both forms
must coincide. This proves the assertion.
\end{proof}

\section{The vectorial case}

For $n=3$, $r=1$, $H_{1}(\mathcal{O})$ may be identified with $h(\mathcal{O})$
by means of the correspondence $u_{1}e^{1}+u_{2}e^{2}+u_{3}e^{3}%
\longleftrightarrow\left(  u_{1},u_{2},u_{3}\right)  $. Also, $P_{1}%
(\mathcal{O})$ may be identified with the space $p(\mathcal{O})$ of
holomorphic vector fields by through the correspondence $w_{0}+w_{2}%
\longleftrightarrow\left(  f,\bar{v}\right)  $, where $f=w_{0}$ and $\bar
{v}=\left(  v_{1},v_{2},v_{3}\right)  $ such that $\ast w_{2}=v_{1}e^{1}%
+v_{2}e^{2}+v_{3}e^{3}$. According to these correspondences the subspace
$\mathcal{N}\cap P_{1}(\mathcal{O})$ may be identified with the subspace of
$p(\mathcal{O})$ defined by
\[
\mathcal{N}^{\prime}:=\left\{  \left(  0,%
\grad
h\right)  :\Delta h=0\text{ \ in }\mathcal{O}\right\}  .
\]
Now it is clear that
\[
\bar{P}_{1,0}(E\smallsetminus K)\cong\bar{p}_{0}(E\smallsetminus
K)=p_{0}(E\smallsetminus K)\diagup\mathcal{N}^{\prime}.
\]
Similarly, $\bar{p}(K)$ may be identified with $\bar{P}_{1}(K)$.

Thus, considering the appropriate topologies we obtain:

\begin{theorem}
Let $K$ be a compact set in $%
\mathbb{R}
^{3}$. Then

1) $\Lambda\in\left(  h(K)\right)  ^{\ast}$ if and only if there exists
$\left(  f,\bar{v}\right)  \in p_{0}(E\smallsetminus K)$ such that
\[
\Lambda\left(  \rho_{\mathcal{O}}(\bar{u})\right)  =-\frac{1}{4\pi}%
\int_{\partial K_{1}}\left\langle \bar{v},\bar{N}\times\bar{u}\right\rangle
dS+\frac{1}{4\pi}\int_{\partial K_{1}}\left\langle \bar{N},f\bar
{u}\right\rangle dS
\]
for every open neiborhood $\mathcal{O}$ of $K$, any vector field $\bar{u}\in
h(\mathcal{O})$ and any regular compact $K_{1}$ satisfying $K\subset
\mathring{K}_{1}\subset K_{1}\subset\mathcal{O}$. Moreover, for every
$\Lambda\in\left(  H_{r}(K)\right)  ^{\ast}$ all the pairs $\left(  f,\bar
{v}\right)  $ with this property belongs to the same class in $\bar{p}%
_{0}(E\smallsetminus K)$.

2) $\Lambda\in\left(  \bar{p}(K)\right)  ^{\ast}$ if and only if there exists
$\bar{u}\in h_{0}(E\smallsetminus K)$ such that%
\[
\Lambda\left(  \rho_{\mathcal{O}}\left(  f,\bar{v}\right)  \right)  =-\frac
{1}{4\pi}\int_{\partial K_{1}}\left\langle \bar{v},\bar{N}\times\bar
{u}\right\rangle dS+\frac{1}{4\pi}\int_{\partial K_{1}}\left\langle \bar
{N},f\bar{u}\right\rangle dS
\]
for every open neiborhood $\mathcal{O}$ of $K$, any pair $\left(  f,\bar
{v}\right)  \in p(\mathcal{O})$ and any regular compact $K_{1}$ satisfying
$K\subset\mathring{K}_{1}\subset K_{1}\subset\mathcal{O}$. Moreover, for a
given $\Lambda\in\left(  \bar{p}(K)\right)  ^{\ast}$ there exists a unique
$\bar{u}\in h_{0}(E\smallsetminus K)$ with this property.
\end{theorem}

Finally, Corollary \ref{cor01} may be reformulated as follows.

\begin{corollary}
If $\bar{u}\in h_{0}(E\smallsetminus K)$ then there exists a scalar function
$f$ and a vector field $\bar{v}$, vanishing both at infinity, such that
$\Delta f=0$, $\Delta\bar{v}=\bar{0}$, $%
\diver
\bar{v}=0$ and
\[
\bar{u}=%
\grad
f+%
\rotor
\bar{v}\text{ \ \ \ in }E\smallsetminus K.
\]

\end{corollary}
\section*{Acknowledgement} We would like to express our gratitude to Victor Petrovich Havin for communicating us
the problems solved here.

\end{document}